\documentclass[12pt,a4paper]{amsart}
\usepackage{amssymb,latexsym}

\theoremstyle{plain}
\newtheorem{lemma}{Lemma}[section]
\newtheorem{theorem}[lemma]{Theorem}

\theoremstyle{definition}
\newtheorem{definition}[lemma]{Definition}
\theoremstyle{remark}
\newtheorem{remark}[lemma]{Remark}
\numberwithin{equation}{section}

\newcommand{\R}{\mathbb{R}}
\renewcommand{\L}{\mathbb{L}}
\newcommand{\C}{\mathbb{C}}
\newcommand{\Ass}{\mathbf{Ass}}
\newcommand{\Con}{\mathbf{Con}}
\newcommand{\Sch}{\mathbf{S}}

\newcommand{\reg}{\text{reg}}
\newcommand{\Ker}{\mathbf{Ker}}

\newcommand{\GL}{\mathbf{GL}}

\newcommand{\E}{\mathbf{E}}
\newcommand{\N}{\mathbf{N}}
\newcommand{\EE}{\mathcal{E}}

\renewcommand{\phi}{\varphi}
\newcommand{\pr}{\mathbf{pr}}

\begin{document}
\title {On Curvatures of Sections of Tensor Bundles.}
\author {P. I. Katsylo}
\address {Independent University of Moscow,
Bolshoi Vlasievskii, 11, 117463, Moscow}
\email {katsylo@katsylo.mccme.ru}
\begin{abstract}
We consider natural differential operations acting on
sections of tensor vector bundles. Arising problems
can be reformulated as invariant theoretical problems
(the IT-reduction). We give examples of usage of the
IT-reduction. In particular, on a manifold with
a connection and a Poisson structure we construct
the canonical quantization.
\end{abstract}
\maketitle

\section{Introduction.}\label{intr}

Let $M$ be a differentiable $m$-dimensional manifold.
The objects of our investigation are \it natural algebraic \rm
differential operations in the form
\[
F_M : \Gamma(\Ass_M(V))^{\reg} \rightarrow \Gamma(\Ass_M(U)),
\]
where $\Ass_M(V)$ and $\Ass_M(U)$ are tensor bundles,
$\Gamma(\Ass_M(U))$ is the space of differentiable sections
of the bundle $\Ass_M(U)$, and
$\Gamma(\Ass_M(V))^{\reg}$ is the space of differentiable
\it nondegenerated \rm sections of the bundle $\Ass_M(V)$.

A tensor bundle $\Ass_M(V)$ is a bundle associated
to the frame bundle on $M$ \cite{Be}, \cite{No}. It
corresponds to some linear representation
$\rho : \GL(m) \rightarrow \GL(V)$. If
\[
\phi : M \rightarrow N
\]
is a diffeomorphism of $M$ onto an open submanifold in $N$, then
\[
\phi^*(\Ass_N(V)) \ = \ \Ass_M(V).
\]
The typical examples of tensor vector bundles are tangent
bundle $T(M)$, cotangent bundle $T^*(M)$, $k$-th exterior
power of cotangent bundle $\wedge^k(M)$ and so on.
The nondegeneracy condition is defined by some
$r(v) \in \R[V]^{\GL(m)}$. A section $s$ of the bundle
$\Ass_M(V)$ is called \it nondegenerated \rm iff for any
chart $Y$ we have
\[
  s_Y(y) \in V^\reg \ \overset{\text{def}}{=}
  \ \{ v \in V \ | \ r(v) \neq 0 \} \subset V
\]
for all $y \in Y$, where $s_Y : Y \rightarrow V$ is
a presentation of $s$ at the chart $Y$.
For example, if $r(v) \equiv \text{const} \neq 0$, then
every section is nondegenerated.

A differential operation $F_M$ is called \it natural \rm
iff for any diffeomorphism
\[
\phi : M \rightarrow N
\]
of $m$-dimensional manifold $M$ onto an open submanifold
of $N$ we have
\[
F_M(\phi^*(\eta)) \ = \ \phi^*(F_N(\eta))
\]
for all $\eta \in \Gamma(\Ass_N(V))^{\reg}$, where
\[
\phi^* : \Gamma(\Ass_N(\cdot)) \rightarrow \Gamma(\Ass_M(\cdot))
\]
is an induced by $\phi$ mapping. This condition is equivalent
to the following one. The presentation of $F_M$ in bundle
charts of $\Ass_M(V)$ and $\Ass_M(U)$
is written by some universal formulas. These formulas
depend on $m$ and on representations $\GL(m) : V, U$.
The condition that $F_M$ is \it algebraic \rm means
that these universal formulas are algebraic. So
for any chart $Y \subset M$ we have: $F_M (s_Y)$ is
a polynomial function in $r(s_Y)^{-1}$ and
\[
\left\{ \frac{\partial^{|i|} s_Y(y)}
{\partial y_1^{i_1} \ldots \partial y_m^{i_m}}
\right\}_{0 \le |i| \le k}.
\]
We consider in this article only algebraic differential
operations. We omit the adjective "algebraic" in the
sequel.

If $F_M$ is a universal differential operation and $s$
is a section of the corresponding tensor bundle, then
$F_M(s)$ is called \it a curvature \rm of that section.

The classical example of a natural differential operation
is the exterior derivative
\[
  d : \Gamma(\wedge^n(M)) \rightarrow \Gamma(\wedge^{n+1}(M)).
\]
$d$ is a linear differential operator of order 1. There is no
other natural linear differential operations of order 1
(the Schouten theorem). Here are some other classical
curvatures and natural differential operations: the commutator
of two vector fields, the curvature of a Riemannian metric,
the Laplacian on a Riemannian manifold \cite{Be}.

In this article we propose the IT-reduction method.
This method gives a reduction of problems about
natural differential operations to the corresponding
invariant theoretical problems. We show how to use the
IT-reduction by considering several known problems.
By the IT-reduction we obtain in \S \ref{metr} a description
of curvatures of a Riemannian metric. (This description is
known in the classical differential geometry \cite{Ep}).
In \S \ref{onriem} we describe differential operations
on a Riemannian manifold.
In \S \ref{conn} we describe curvatures of a connection.
In fact, a connection is not a section of any tensor
bundle. But one can apply to connections the modified 
IT-reduction method. In \S \ref{withconn} we describe
differential operations on a manifold with a connection.
As a corollary of this description we obtain in \S \ref{q}
that for any manifold with a connection $\theta$ and a
Poisson structure $\omega$ there exists the canonical
quantization $\star(\theta, \omega)$. The word "canonical"
means that if $N$ is a manifold with a connection $\theta$
and a Poisson structure $\omega$ and
\[
\phi : M \rightarrow N
\]
is a diffeomorphism of a manifold $M$ onto an open
submanifold of $N$, then
\[
  \star(\phi^\ast(\theta), \phi^\ast(\omega)) \ =
  \ \phi^*(\star(\theta, \omega)).
\]

Finally, we propose the following Conjecture.

\noindent \bf Conjecture. \it
Let $M$ be a manifold and $\Gamma(\Ass_M(\wedge^2 V))^P$
be a set of all Poisson structures on $M$. Then there are no
natural differential operations
\[
  \Xi_k \ : \ \Gamma(\Ass_M(\wedge^2 V))^P \times
  C^\infty(M) \times C^\infty(M)
  \rightarrow C^\infty(M), \qquad
  k = 1,2, \ldots
\]
such that the operation
\begin{equation*}\begin{split}
  \star \ : \ C^\infty(M) \times C^\infty(M)
  \rightarrow C^\infty(M)[[\hbar]], \\
  (f,g) \mapsto f g + \Xi_1(\omega, f, g) \hbar +
  \Xi_2(\omega, f, g) \hbar^2 + \ldots .
\end{split}\end{equation*}
is a quantization of the Poisson structure $\omega$. \rm

Roughly speaking, this conjecture claims that for a
canonical quantization of a Poisson structure on a manifold
one needs some additional structure on that manifold.

\section{Invariant theoretical description of
natural differential operations.}\label{descr}

Let $x^1, \ldots , x^m$ and $z^1, \ldots , z^m$ be two copies
of the standard basis of the space $\R^{m*}$. We consider
$x^1 , \ldots , x^m$ as coordinate functions in $\R^m$.
Consider the linear space $\EE_k$ of $k$-jets of germs of
differentiable functions at $0 \in \R^m$.
We identify the space $\EE_k$ and the space of polynomials
in the variables $z^1, \ldots , z^m$ of degree $\le k$.
Consider the group $\GL(m)_{k+1}$ of $(k+1)$-jets
of germs of diffeomorphisms of a neighborhood of $0$ in $\R^m$
onto a neighborhood of $0$ in $\R^m$. By $\GL(m)_\infty$ we
denote the group of germs of diffeomorphisms of a neighborhood
of $0$ in $\R^m$ onto a neighborhood of $0$ in $\R^m$.

We have the canonical group homomorphism
\[
  \phi_{k+1} : \GL(m)_{k+1} \rightarrow \GL(m),
\]
where $\phi_{k+1}(g)$ is the Jacobi matrix of $g$ at $0$.
Set
\[
  \N(m)_{k+1} \ =
  \ \{ g \in \GL(m)_{k+1} \ | \ \phi_{k+1}(g) = \E_m \},
\]
where $\E_m$ is the identity matrix of size $m \times m$.
The group $\GL(m)_{k+1}$ is a real linear algebraic group,
$\N(m)_{k+1}$ is its unipotent radical, and
$\GL(m)_{k+1} / \N(m)_{k+1} \simeq \GL(m)$.

Let $\rho : \GL(m) \rightarrow \GL(V)$ be a linear representation.
Consider the space $\EE_k \otimes V$ as a linear space of $k$-jets
of germs of differentiable mappings of a neighborhood of
$0$ in $\R^m$ to $V$. By $\EE_\infty \otimes V$ we denote
the linear space of germs of differentiable
mappings of a neighborhood of $0$ in $\R^m$ to $V$.
The group $\GL(m)_\infty$ acts canonically on the space
$\EE_\infty \otimes V$:
\[
  \GL(m)_\infty \ : \ \EE_\infty \otimes V, \qquad
  (g \cdot v)(x) = \rho(J_g(g^{-1}(x))(v(g^{-1}(x))),
\]
where $J_g(x)$ is the Jacobi matrix of $g$ at $x$.
This action corresponds to the transition rule for local
presentations of a section of the bundle $\Ass_M(V)$ in
bundle charts. It defines canonically the action
\begin{equation}\label{E:jact}\begin{split}
  \GL(m)_{k+1} \ : \ \EE_k \otimes V, \qquad
  (g \cdot v)(z) = \\
  \{ k \text{-jet of the mapping}
  \ \ x \mapsto \rho(J_g(g^{-1}(x))(v(g^{-1}(x))) \},
\end{split}\end{equation}
where $g = g(z) \in \GL(m)_{k+1}$,
\ $v = v(z) \in \EE_k \otimes V$.

Note that for $a \ge b$ we have the canonical group
homomorphism
\[
  \GL(m)_a \ \rightarrow \ \GL(m)_b.
\]
The representation \eqref{E:jact} and this homomorphism define
the representation
\[
\GL(m)_a : \EE_b \otimes V
\]
for all $a > b$.

We have the following $\GL(m)_{k+1}$-equivariant linear mapping
\[
  \varepsilon_k : \EE_k \otimes V \rightarrow V,
  \qquad v(z) \mapsto v(0).
\]

Fix $r(v) \in \R[V]^{\GL(m)}$ and consider
$\GL(m)$-invariant subset
\[
  V^\reg \ = \ \{ v \in V \ | \ r(v) \neq 0 \} \subset V.
\]
This subset defines the natural algebraic nondegeneracy condition
for sections of the tensor bundle $\Ass_M(V)$. Set
\[
  (\EE_k \otimes V)^{\reg} = \varepsilon_k^{-1}(V^{\reg}).
\]
\begin{theorem}\label{T:main}
  Let $\sigma : \GL(m) \rightarrow \GL(U)$ be a linear
  representation. Then there is the canonical
  1-1 correspondence between $\GL(m)_{k+1}$-equivariant
  morphisms in the form
  \[
    \alpha : \ (\EE_k \otimes V)^\reg \ \rightarrow
    \ U = \EE_0 \otimes U
  \]
  and the set of natural differential operations of
  $k$-th order in the form
  \[
    F_{\alpha, M} : \Gamma(\Ass_M(V))^{\reg}
    \ \rightarrow \ \Gamma(\Ass_M(U)).
  \]
\end{theorem}
\begin{proof}
  First, we show that any $\GL(m)_{k+1}$-equivariant morphism
  defines canonically a natural differential operation.
  
  Let \ $\alpha : \EE_k \otimes V \rightarrow U$ \ be
  \ $\GL(m)_{k+1}$-equivariant morphism. Consider local
  coordinates $x = (x^1, \ldots ,x^m)$ in some
  neighborhood of a point $p \in M$.
  We assume that the local coordinates defines 1-1
  mapping of an open subset $p \in X \subset M$ and some
  neighborhood $\widetilde{X}$ of $0$ in $\R^m$.
  Define the natural differential operation $F_{\alpha, M}$
  in the following way. Suppose that
  $s_X : \widetilde{X} \rightarrow V$
  is a presentation of a section $s$ of the bundle
  $\Ass_M(V)$ in the chart $X$. We define the section
  $F_{\alpha, M}(s)$ of the bundle $\Ass_M(U)$ in the chart
  $X$ by the following formula:
  \begin{gather}
    (F_{\alpha, M}(s))_X : \widetilde{X}
    \ \rightarrow \ U, \notag \\
    x \ \mapsto \ \alpha\left(\sum_{0 \le |i| \le k}
    \frac{|i|!}{i_1! \ldots i_m!}\frac{\partial^{|i|}s_X(x)}
    {\partial (x^1)^{i_1} \ldots \partial (x^m)^{i_m}}
    (z^1)^{i_1} \ldots (z^m)^{i_m}\right). \notag
  \end{gather}
  Let us check that the section $F_{\alpha, M}(s)$ is
  well-defined. Let $y = (y^1, \ldots ,y^m)$ be some other
  local coordinates in some neighborhood $Y$ of the point
  $p$. The transition formulas from the local coordinates $x$
  to the local coordinates $y$ define a diffeomorphism $g$ of a
  neighborhood of $0$ in $\R^m$ onto a neighborhood of $0$ in
  $\R^m$. Let $g_{k+1}$ be $(k+1)$-jet of the diffeomorphism $g$
  at $0 \in \R^m$. Suppose that
  $s_Y : \widetilde{Y} \rightarrow V$ is a
  presentation of the section $s$ in the chart $Y$. By the
  transition rule we have
  \begin{equation}\label{E:transm}
    s_Y(g(x)) = \rho (J(x))(s_X(x)),
  \end{equation}
  where $J(x)$ is the Jacobi matrix of the mapping $g$.
  We have to check that
  \[
    (F_{\alpha, M}(s))_Y(0) =
    \sigma (J(0)) ((F_{\alpha, M}(s))_X(0)).
  \]
  Consider $k$-jets of the left and the right sides of the
  equation \eqref{E:transm}. We obtain
  \begin{gather*}
    \sum_{0 \le |i| \le k}\frac{|i|!}{i_1! \ldots i_m!}
    \left.\frac{\partial^{|i|}s_Y(y)}
    {\partial (y^1)^{i_1} \ldots \partial (y^m)^{i_m}}
    \right|_{y=0}
    (z^1)^{i_1} \ldots (z^m)^{i_m} = \\
    g_{k+1} \cdot \left(\sum_{0 \le |i| \le k}
    \frac{|i|!}{i_1! \ldots i_m!}
    \left.\frac{\partial^{|i|}s_X(x)}
    {\partial (x^1)^{i_1} \ldots \partial (x^m)^{i_m}}
    \right|_{x=0}
    (z^1)^{i_1} \ldots (z^m)^{i_m} \right).
  \end{gather*}
  By the $\GL(m)_{k+1}$-equivalence of the morphism
  $\alpha$ we obtain
  \begin{gather*}
  (F_{\alpha, M}(s))_Y(0) = \\
    \alpha\left(g_{k+1} \cdot \left(\sum_{0 \le |i| \le k}
    \frac{|i|!}{i_1! \ldots i_m!}
    \left.\frac{\partial^{|i|}s_X(x)}
    {\partial (x^1)^{i_1} \ldots \partial (x^m)^{i_m}}
    \right|_{x=0}
    (z^1)^{i_1} \ldots (z^m)^{i_m} \right)\right) = \\
    \sigma (J(0))
    \left(\alpha\left(\sum_{0 \le |i| \le k}
    \frac{|i|!}{i_1! \ldots i_m!}
    \left.\frac{\partial^{|i|}s_X(x)}
    {\partial (x^1)^{i_1} \ldots \partial (x^m)^{i_m}}
    \right|_{x=0}
    (z^1)^{i_1} \ldots (z^m)^{i_m} \right)\right) = \\
    \sigma (J(0))((F_{\alpha, M}(s))_X(0)).
  \end{gather*}
  
  Conversely, by the construction above any natural differential
  operation defines canonically the $\GL(m)_{k+1}$-morphism.
\end{proof}

Theorem \ref{T:main} reduces the problem of description of
natural differential operations of order $k$ acting on
$\Gamma(\Ass_M(V))^\reg$ to an invariant theoretical problem.
Namely, consider the regular action
\[
  \N(m)_{k+1} : (\EE_k \otimes V)^\reg
\]
and the corresponding algebra of invariants
\[
  \R[(\EE_k \otimes V)^\reg]^{\N(m)_{k+1}}.
\]
We have the canonical representation
\[
  \GL(m)_{k+1}/\N(m)_{k+1} \simeq \GL(m) \ :
  \ \R[(\EE_k \otimes V)^\reg]^{\N(m)_{k+1}}.
\]
Every $\GL(m)$-embedding
\begin{equation}\label{incl}
  A_U \ : \ U \ \hookrightarrow
  \ \R[(\EE_k \otimes V)^\reg]^{\N(m)_{k+1}}
\end{equation}
of finite dimensional $\GL(m)$-module $U$ defines canonically
$\GL(m)_{k+1}$-morphism
\begin{equation}\label{morf}
    \alpha_U : \ (\EE_k \otimes V)^{\reg} \ \rightarrow
    \ U^\ast = U^\ast \otimes \EE_0, \qquad
    \alpha_U(v(z))(u) = A_U(u)(v(z)).
\end{equation}
Conversely, every $\GL(m)_{k+1}$-morphism
\eqref{morf} corresponds to some embedding \eqref{incl}.

\section{Curvatures of the Riemannian metrics.}\label{metr}

Consider $m$-dimensional manifold $M$.
Let $e_1 , \ldots , e_m$ be the standard basis of $\R^m$ and
$z^1 , \ldots , z^m$ and $u^1 , \ldots , u^m$ are two copies
of the dual basis of the dual space $\R^{m\ast}$.
We identify the space $\EE_k$ and the space of polynomials
in the variables $z^1 , \ldots , z^m$ of degree $\le k$.
Riemanian metrics are nondegenerated sections of the bundle
$\Ass_M(S^2 \R^{ m \ast})$. Fix $k$ and consider the action
\begin{equation}\label{eqn}
  \N(m)_{k+1} \ : \ (\EE_k \otimes S^2 \R^{m\ast})^\reg. 
\end{equation}
In this article we need a simplified concept of
\it the Seshadri section \rm (see \cite{Se}).
\begin{definition}
  Let $G$ be a linear algebraic group,
  $G:X$ be a regular action on an affine variety.
  A closed subvariety $Y \subset X$ is called the
  \it nice Seshadri section, \rm if
  \begin{itemize}
  \item $G Y = X$,
  \item every $G$-orbit intersects transversally $Y$ 
    at one point.
  \end{itemize}
\end{definition}
\noindent Suppose that $Y \subset X$ is a nice Seshadri section
and $V$ is a vector space. Consider $V$ as a trivial
$G$-module. Then we have the canonical 1-1 correspondence
\begin{gather*}
  \{ \text{set of regular mappings of} \ Y
  \ \text{to} \ V \} \rightarrow \\
  \{ \text{set of regular} \ G-\text{mappings of}
  \ X \ \text{to} \ V \}, \\
  \{ \hat{\xi} : Y \rightarrow V \} \longrightarrow
  \{ \xi : X \rightarrow V,
  \ \ \xi(x) = \hat{\xi}((G \cdot x) \cap Y) \}.
\end{gather*}

Our next purpose is to construct the nice Seshadri section
for the action \eqref{eqn}.

Consider the group $\N(m)_{k+1}$ as an affine
variety. Then it is isomorphic to a linear space.
Namely, we have the following isomorphism
\begin{gather*}
  \eta : (S^2 \R^{m\ast} \otimes \R^m) \times
  (S^3 \R^{m\ast} \otimes \R^m) \times \ldots \times
  (S^{k+1} \R^{m\ast} \otimes \R^m) \ \rightarrow \ \N(m)_{k+1}, \\
  (g_2, g_3, \ldots , g_{k+1}) \mapsto
  \eta(g_2, g_3, \ldots , g_{k+1}) =
  (\E + g_{k+1}) \cdot \ldots \cdot (\E + g_3) \cdot (\E + g_2),
\end{gather*}
where
$g_n = g_{ni}(z) \otimes e_i \in S^n \R^{m\ast} \otimes \R^m$,
\ $g_{ni}(z) \in S^n \R^{m\ast}$, \ $\E + g_n$ \ is
$k$-jet at $0$ of the mapping
\[
  \R^m \rightarrow \R^m, \qquad c_i e_i \mapsto
  (c_i + g_{ni}(c_1 , \ldots , c_m)) e_i .
\]
We use the following identification
\begin{gather*}
  (\EE_k \otimes S^2 \R^{m\ast})^\reg \simeq
  (S^2 \R^{m\ast})^\text{reg} \times
  (\R^{m\ast} \otimes S^2 \R^{m\ast}) \times \ldots
  \times (S^k \R^{m\ast} \otimes S^2 \R^{m\ast}), \\
  h \ = \ h_0 + h_1 + \ldots + h_k \ \sim
  \ (h_0, h_1, \ldots , h_k),
\end{gather*}
where $h_0 \in (S^2 \R^{m\ast})^\text{reg}$,
\ $h_n = h_n(z,u)$ is a bihomogeneous polynomial,
$\deg_z(h_n) = n$, \ $\deg_u(h_n) = 2$, i.e.,
$h_n \in S^n \R^{m\ast} \otimes S^2 \R^{m\ast}$.

Suppose
\[
  h = (h_0, h_1, \ldots , h_k) \in
  (\EE_k \otimes S^2 \R^{m\ast})^\reg.
\]
For $g = \eta(g_2, g_3, \ldots , g_{k+1}) \in \N(m)_{k+1}$
we have
\[
  g \cdot h = ((g \cdot h)_0, (g \cdot h)_1,
  \ldots ,(g \cdot h)_k),
\]
where $(g \cdot h)_n = \mu_n (\{ g_i, h_j \})$.
From \eqref{E:jact} it is not difficult to obtain that
\[
  \mu_n (\{ g_i, h_j \}) =  \mu'_n(g_{n+1}, h_0) +
  \mu''_n(g_2, \ldots , g_n, h),
\]
where
\begin{equation}\label{E:quasi}
  \mu'_n(g_{n+1}, h_0) \ =
  \ -\frac{\partial^2 g_{n+1}}{\partial e_i \partial z_j}
  \otimes \frac{\partial h_0}{\partial u^i} u^j.
\end{equation}

For $n \ge 2$ define $\GL(m)$-submodule
\[
   L_n \ = \ \Ker(\delta_n) \ \subset
   \ S^n \R^{m\ast} \otimes S^2 \R^{m\ast},
\]
where
\begin{gather*}
  \delta_n \ : \ S^n \R^{m\ast} \otimes S^2 \R^{m\ast}
  \ \rightarrow \ S^{n+1} \R^{m\ast} \otimes \R^{m\ast}, \\
  \delta_n (f(z) \otimes q(u)) =
  f(z) z^i \otimes \frac{\partial q(u)}{\partial u^i}.
\end{gather*}
The representation $\GL(m) : L_n$ is irreducible.
It is isomorphic to the representation
$\GL(m) : \Sch_{(n,2)}(\R^{m\ast})$, where $\Sch_\lambda$ is the
Schur functor corresponding to the partition $\lambda$
(see \cite{Fu}).
Consider the subvariety
\[
  \L_k \overset{\text{def}}{ = }
  (S^2 \R^{m\ast})^\reg \times \{ 0 \}
   \times
  L_2 \times L_3 \times \ldots \times L_k
  \subset (\EE_k \otimes S^2 \R^{m\ast})^\reg.
\]
\begin{lemma}\label{sesh}
  Every $\N(m)_{k+1}$-orbit intersects transversally
  $\L_k$ at one point. In other words, $\L_k$ is
  a nice Seshadri section for the action \eqref{eqn}.
\end{lemma}
\begin{proof}
  Suppose $h \in (\EE_k \otimes S^2 \R^{m\ast})^\reg$.
  For $g = \eta(g_2, g_3, \ldots , g_{k+1}) \in \N(m)_{k+1}$
  the condition $g \cdot h \in \L_k$ is equivalent to
  the equations ($E_1$) - ($E_k$), where
  \begin{enumerate}
  \item[($E_1$)] $((\E + g_2) \cdot h)_1 = 0$,
  \item[($E_n$)] $((\E + g_{n+1}) \cdot \ldots \cdot (\E + g_3)
    \cdot (\E + g_2) \cdot h))_n \in L_n$,
    where $2 \le n \le k$.
  \end{enumerate}
  We claim that the equations ($E_1$) - ($E_k$) for $g$ have
  a unique solution. More precisely:
  
  \it $(\ast)$ One can find sequentially the elements
  $g_2, \ldots , g_{k+1}$ in a unique way from the equations
  $(E_1), \ldots , (E_k)$ accordingly. Moreover, the equation
  ($E_n$) for $g_{n+1}$ with fixed (before defined)
  $g_2, \ldots , g_n$ is a linear equation that has a unique
  solution. \rm
  
  Let us prove $(\ast)$. First, by the nondegeneracy of $h_0$
  and the $\GL(m)$-equivalence we can assume that
  $h_0 = (u^1)^2 + \ldots + (u^m)^2$.
  
  Consider the equation $(E_1)$ for $g_2$. By \eqref{E:quasi} we
  can rewrite it in the following way:
  \[
    -2 \frac{\partial^2 g_2}{\partial e_i \partial z^j}
  \otimes  u^i  u^j + h_1 \ = \ 0.
  \]
  It is easy to see that this equation for $g_2$ is a linear
  equation having a unique solution.
  
  Suppose that we find $g_2, \ldots , g_n$ from the equations
  $(E_1), \ldots , (E_{n-1})$. Consider the equation $(E_n)$
  for $g_{n+1}$. By \eqref{E:quasi} and the definition of $L_n$
  we can rewrite the equation $(E_n)$ in the following way:
  \begin{equation}\label{eq1}
    \delta_n\left( -2 \frac{\partial^2 g_{n+1}}
    {\partial e_i \partial z^j}
  \otimes  u^i  u^j + h_n'\right) \ = \ 0,
  \end{equation}
  where $h_n' = h_n'(h, g_2, \ldots , g_n) \in
  S^n \R^{m\ast} \otimes S^2 \R^{m\ast}$. Using the definition
  of $\delta_n$ we simplify the equation
  \eqref{eq1} to the following one.
  \begin{equation}\label{eq2}
    -2 \left(\frac{\partial^2 g_{n+1}}{\partial e_i \partial z^j}+
      \frac{\partial^2 g_{n+1}}{\partial e_j \partial z^i}\right)
    z^i \otimes u^j + \delta_n(h_n') \ = \ 0.
  \end{equation}
  It is not difficult to check that the equation \eqref{eq2}
  for $g_{n+1}$ is a linear equation having a unique solution.
\end{proof}

Consider the mapping
\[
  \widetilde{\alpha}_k : (\EE_k \otimes S^2 \R^{m\ast})^\reg
  \ \rightarrow \ \L_k, \qquad
  h \mapsto (\N(m)_{k+1} \cdot h) \cap \L_k.
\]
We have the natural action of the group $\GL(m)$ on $\L_k$.
This action defines canonically the action
$\GL(m)_{k+1} : \L_k$ such that the subgroup $\N(m)_{k+1}$
acts on $\L_k$ trivially.

For $n = 2,3, \ldots ,k$ let
\[
  \pr_n : \L_k \rightarrow L_n
\]
be the canonical projections. Then
\[
  \alpha_n \overset{\text{def}}{ = }
  \pr_n \circ \widetilde{\alpha}_k \ :
  \ (\EE_k \otimes S^2 \R^{m\ast})^\reg
  \rightarrow L_n
\]
is $\GL(m)_{k+1}$-equivariant morphism.
By Theorem \ref{T:main} \ $\alpha_n$ corresponds to
a natural differential operation
\[
  A_n \ : \ \Gamma(\Ass_M(S^2 \R^{m\ast}))^\reg \ \rightarrow
  \ \Gamma(\Ass_M(L_n)).
\]
Lemma \ref{sesh} implies the following statement.
\begin{theorem}
  Let $U$ be $\GL(m)$-module and
  \[
    F \ : \ \Gamma(\Ass_M(S^2 \R^{m\ast}))^\reg \ \rightarrow
    \ \Gamma(\Ass_M(U))
  \]
  be a natural differential operation of order $k$.
  Then
  \[
    F(h) = \widetilde{F}
    (\det(h)^{-1}, h, A_2(h), \ldots , A_k(h)),
  \]
  where $\widetilde{F}$ corresponds to some polynomial
  $\GL(m)$-mapping
  \[
    \widetilde{f} \ : \ \R \times S^2 \R^{m\ast} \times
    L_2 \times \ldots \times L_k \ \rightarrow \ U.
  \]
\end{theorem}
\begin{remark}
  The classical description of the operation $A_n$
  is the following: $A_n(h)|_p$ is the homogeneous of
  degree $n$ summand of the Taylor series of the metric
  $h$ at the normal local coordinates with center $p$
  (see \cite{Ep}).
\end{remark}

\section{Natural differential operations on the Riemannian
manifolds.}\label{onriem}

\begin{theorem}\label{mainonriem}
  Let $\rho : \GL(m) \rightarrow \GL(V)$
  and $\sigma : \GL(m) \rightarrow \GL(U)$ be linear
  representations. Then there is the canonical 1-1 correspondence
  between $\GL(m)$-equivariant differential
  operations of bounded order in the form
  \begin{equation}\label{opriem}
    \xi : \L_k  \times C^\infty(\R^m, V)
    \rightarrow C^\infty(\R^m, U),
  \end{equation}
  where $k \in \N$ and the set of natural differential
  operations of bounded order in the form
  \[
    \Xi : \Gamma(\Ass_M(V))
    \ \rightarrow \ \Gamma(\Ass_M(U))
  \]
  on a Riemannian manifold.
\end{theorem}
\begin{proof}
  Let $z^1 , \ldots , z^m$ be the standard basis of the space
  $\R^{m\ast}$. We identify $\EE_k$ and the space of polynomials
  in the variables $z^1 , \ldots , z^m$ of degree $\le k$.
  
  Suppose we have $\GL(m)$-equivariant differential
  operation \eqref{opriem} of an order $\le k$.
  Consider the regular action
  \[
    \N(m)_{k+1} : (\EE_k \otimes S^2(\R^{m\ast}))^\reg
    \times (\EE_k \otimes V).
  \]
  Note that the subvariety $\L_k \times (\EE_k \otimes V)$
  is the nice Seshadri section for this action. Define
  $\GL(m)$-mapping
  \[
    \hat{\xi} \ : \ \L_k \times (\EE_k \otimes V)
    \rightarrow U, \qquad
    (h, v(z)) \mapsto \xi(h, v(z))(0).
  \]
  Extend the mapping $\hat{\xi}$ to the $\GL(m)_{k+1}$-mapping
  \[
    \xi \ : \ (\EE_k \otimes S^2(\R^{m\ast}))^\reg
    \times (\EE_k \otimes V)
    \ \rightarrow \ U.
  \]
  From Theorem \ref{T:main} it follows that $\xi$ defines
  the natural differential operation
  \[
    \tilde{\Xi} \ : \ \Gamma(\Ass_M(S^2 \R^{m\ast}))^\reg
    \times \Gamma(\Ass_M(V))
    \ \rightarrow \ \Gamma(\Ass_M(U)).
  \]
  Now suppose that $(M, h)$ is a Riemannian manifold.
  Define the operation
  \[
    \Xi \ : \Gamma(\Ass_M(V))
    \ \rightarrow \ \Gamma(\Ass_M(U)), \qquad
    s \mapsto \tilde{\Xi}(h, s).
  \]

  Conversely, by the construction above any natural differential
  operation on a Riemannian manifold defines canonically the
  $\GL(m)$-equivariant differential operation \eqref{opriem}
  of bounded order.
\end{proof}
For example, the Laplacian
\[
  \Delta \ : \ \Gamma(\Ass_M(V)) \ \rightarrow
  \Gamma(\Ass_M(V))
\]
corresponds to the following differential operation
\begin{gather*}
  \delta \ : \ S^2(\R^{m\ast})^\reg \times
  C^\infty(\R^m, V) \ \rightarrow \ C^\infty(\R^m, V), \\
  (h, v(x)) \mapsto (h^{-1})^{ij}
  \frac{\partial^2 v(x)}{\partial x^i \partial x^j},
\end{gather*}
where $h^{-1} \in S^2(\R^m)^\reg$ is the dual to $h$
quadratic form.
\begin{remark}
  From the proof of Theorem \ref{mainonriem} we obtain the following
  description of the corresponding to $\xi$ differential operation
  $\Xi$. Let $M$ be a Riemannian manifold, $s$ be a section of
  $\Ass_M(V)$, and $p \in M$. Take normal local coordinates with
  center $p$ and calculate $\Xi(s)(p)$ in that local coordinates
  by the formula \eqref{opriem}.
\end{remark}

\section{Curvatures of connections.}\label{conn}

Consider $m$-dimensional manifold $M$. Let $\Con(M)$ be
the set of all connections on $M$. In \S \ref{intr} we
define natural differential operations acting on sections
of tensor bundles. Analogously one can define natural
differential operations acting on connections, pairs
$(\theta, \eta)$, where $\theta$ is a connection and
$\eta$ is a section of a tensor bundle.

By definition, a curvature of a connection
$\theta \in \Con(M)$ is $R(\theta)$, where
\[
  R \ : \ \Con(M) \rightarrow \Gamma(\Ass_M(U))
\]
is some natural differential operation.

In this section we describe curvatures of a connection on
a manifold.

Let $e_1 , \ldots , e_m$ be the standard basis of $\R^m$,
$x^1 , \ldots , x^m$, $z^1 , \ldots , z^m$,
$u^1 , \ldots , u^m$, and $v^1 , \ldots , v^m$
be copies of the dual basis of the dual space
$\R^{m\ast}$. We consider $x^1 , \ldots , x^m$ as
coordinate functions in $\R^m$. Define
$\EE_n$, $\GL(m)_n$, $\N(m)_n$ as in \S \ref{descr}.
We identify $\EE_k$ and the space of polynomials
in the variables $z^1 , \ldots , z^m$ of degree $\le k$.

A connection on $0 \in X \subset \R^m$ is a mapping
\[
  \theta \ : \ X \rightarrow
  \R^m \otimes \R^{m\ast} \otimes \R^{m\ast},
  \qquad x \ \mapsto \ \theta(x) =
  \theta_{ij}^l(x) e_l \otimes u^i \otimes v^j.
\]
The mapping $\theta$ defines the connection operator $D_\theta$
in the following way:
\[
  D_\theta(e_i) \ = \ \theta_{ji}^l(x) e_l \otimes u^j.
\]

Let $\Con(m)_\infty$ be the set of germs of connections
at $0 \in \R^m$ and $\Con(m)_k$ be the set of $k$-jets of
germs of connections at $0 \in \R^m$. We use the following
identification
\begin{gather*}
  \Con(m)_k \ \simeq \ (\R^m \otimes \R^{m\ast} \otimes \R^{m\ast})
  \times
  (\R^{m\ast} \otimes \R^m \otimes \R^{m\ast} \otimes \R^{m\ast})
  \times \\
  (S^2(\R^{m\ast}) \otimes \R^m \otimes \R^{m\ast} \otimes \R^{m\ast})
  \times \ldots \times
  (S^k(\R^{m\ast}) \otimes \R^m \otimes \R^{m\ast} \otimes \R^{m\ast}),\\
  \theta = \theta_0 + \theta_1 + \ldots + \theta_k \ \sim
  \ (\theta_0, \theta_1, \ldots , \theta_k),
\end{gather*}
where $\theta_n = \theta_n(z, e, u, v)$ is a polyhomogeneous
polynomial of polydegree $(n,1,1,1)$, i.e., $\theta_n \in
S^n(\R^{m\ast}) \otimes \R^m \otimes \R^{m\ast} \otimes \R^{m\ast}$.

The group $\GL(m)_\infty$ acts canonically on $\Con(m)_\infty$:
\begin{gather*}
  \GL(m)_\infty \ : \ \Con(m)_\infty, \\
  (g \ast \theta)(x) \ = \ \left(J_g(x)_a^l
  \frac{\partial J_g^{-1}(x)_j^a}{\partial x^i} +
  J_g(x)_a^l \theta_{ib}^a(x) J_g^{-1}(x)_j^b \right)
  e_l \otimes u^i \otimes v^j,
\end{gather*}
where $J_g(x)$ is the Jacobi matrix of $g$ at $x$.
This action corresponds to the transition rule for local
presentations of a connection in charts.
It defines canonically the action
\begin{equation}\label{actcon}\begin{split}
  \GL(m)_{k+2} : \Con(m)_k, \qquad
  (g \cdot \theta)(z) = \{ k-\text{jet of the mapping} \\
  x \mapsto
  \left(J_g(x)_a^l
  \frac{\partial J_g^{-1}(x)_j^a}{\partial x^i} +
  J_g(x)_a^l \theta_{ib}^a(x) J_g^{-1}(x)_j^b \right)
  e_l \otimes u^i \otimes v^j \},
\end{split}\end{equation}
where $g = g(z) \in \GL(m)_{k+2}$,
\ $\theta = \theta(z) \in \Con(m)_k$.
The action $\GL(m)_{k+2} : \Con(m)_k$ is an \it affine \rm
action. Consider the restriction of this action to the subgroup
$\N(m)_{k+2} \subset \GL(m)_{k+2}$:
\begin{equation}\label{auconn}
  \N(m)_{k+2} : \Con(m)_k.
\end{equation}
\begin{theorem}\label{T:main2}
  Let $\rho : \GL(m) \rightarrow \GL(V)$ and
  $\sigma : \GL(m) \rightarrow \GL(U)$ be linear representations.
  Then there is the canonical 1-1 correspondence
  between $\GL(m)_{k+2}$-morphisms in the form
  \[
    \alpha : \ \Con(m) \times (\EE_k \otimes V)
    \ \rightarrow \ U = \EE_0 \otimes U
  \]
  and the set of natural differential operations
  of order $k$ in the form
  \[
    F_{\alpha, M} : \Con(M) \times \Gamma(\Ass_M(V))
    \ \rightarrow \ \Gamma(\Ass_M(U)).
  \]
\end{theorem}
The proof of this Theorem is the same as 
of Theorem \ref{T:main}.

Our next purpose is to construct the nice Seshadri section
for the action \eqref{auconn}.

We use the isomorphism $\eta$ from the section \ref{metr}.
Suppose
\[
  \theta = (\theta_0, \theta_1, \ldots , \theta_k)
  \in \Con(m)_k,
\]
where $\theta_n \in
S^n(\R^{m\ast}) \otimes \R^m\otimes \R^{m\ast} \otimes \R^{m\ast}$.
For
$g = \eta(g_2, g_3, \ldots , g_{k+2}) \in \N(m)_{k+2}$
we have
\[
  g \cdot \theta = ((g \cdot \theta)_0, (g \cdot \theta)_1,
  \ldots ,(g \cdot \theta)_k),
\]
where $(g \cdot h)_n = \nu_n (\{ g_i, h_j \})$.
From \eqref{actcon} it is easy to obtain that

\[
  \nu_n(\{ g_i , \theta_j \}) = \nu_n'(g_{n+2}) +
  \nu_n''(g_2, \ldots , g_{n+1}, \theta),
\]
where
\begin{equation}\label{nu}
  \nu_n'(g_{n+2}) \ =
  \ -\frac{\partial^2 g_{n+2}}{\partial z^i \partial z^j}
  \otimes u^i \otimes v^j.
\end{equation}
Define $\GL(m)$-submodule
\[
   C_n \ = \ \Ker(\gamma_n) \ \subset
   \ S^n(\R^{m\ast}) \otimes \R^m \otimes \R^{m\ast}
   \otimes \R^{m\ast},
\]
where
\begin{gather*}
  \gamma_n \ :
  \ S^n(\R^{m\ast}) \otimes \R^m \otimes \R^{m\ast} \otimes \R^{m\ast}
  \ \rightarrow
  \ S^{n+2} \R^{m\ast} \otimes \R^m, \\
  \gamma_n (f(z) \otimes e_l \otimes u^i \otimes v^j)
  = f(z) z^i z^j \otimes e_l.
\end{gather*}
Consider the subvariety
\[
  \C_k \overset{\text{def}}{ = }
  C_0 \times C_1 \times \ldots \times C_k
  \subset \Con(m)_k.
\]
\begin{lemma}\label{sesh2}
  Every $\N(m)_{k+2}$-orbit intersects transversally $\C_k$ 
  at one point. In other words, $\C_k$ is a nice Seshadri
  section for the action \eqref{auconn}.
\end{lemma}
\begin{proof}
  Suppose $\theta \in \Con(m)_k$.
  For $g = \eta(g_2, g_3, \ldots , g_{k+2}) \in \N(m)_{k+2}$
  the condition $g \cdot \theta \in \C_k$ is equivalent
  to equations ($D_0$) - ($D_k$), where
  \begin{enumerate}
  \item[($D_n$)] $((\E + g_{n+2}) \cdot \ldots \cdot (\E + g_3)
    \cdot (\E + g_2) \cdot \theta))_n \in C_n$.
  \end{enumerate}
  It is easy to see that the Lemma is a corollary of the
  following claim.

  \it $(\ast)$ One can find sequentially the elements
  $g_2, \ldots , g_{k+2}$ in a unique way from the equations
  $(D_0), \ldots , (D_k)$ accordingly. Moreover, the equation
  ($D_n$) for $g_{n+2}$ with fixed (before defined)
  $g_2, \ldots , g_{n+1}$ is a linear equation that has a unique
  solution. \rm
  
  Let us prove claim $(\ast)$. By \eqref{nu} and the
  definition of $C_n$ we can rewrite the equation $(D_n)$
  in the following way:
  \[
    \gamma_n\left(-\frac{\partial^2 g_{n+2}}
    {\partial z^i \partial z^j}
    \otimes u^i \otimes v^j  + \theta_n'\right) \ = \ 0,
  \]
  where $\theta_n' = \theta_n'(\theta, g_2, \ldots , g_{n+1})
  \in S^n(\R^{m\ast}) \otimes \R^m \otimes \R^{m\ast}
  \otimes \R^{m\ast}$. Using the definition of $\gamma_n$
  and the Euler theorem about homogeneous functions we get
  \[
    -(n+2)(n+1)g_{n+2} + \gamma_n(\theta_n') \ = \ 0.
  \]
  It is clear that this equation for $g_{n+2}$ with a fixed
  (before defined) $g_2, \ldots , g_{n+1}$ is a
  linear equation having a unique solution.
\end{proof}
Consider the mapping
\[
  \widetilde{\psi}_k : \Con(m)_k
  \ \rightarrow \ \C_k, \qquad
  \theta \mapsto (\N(m)_{k+2} \cdot \theta) \cap \C_k.
\]
The natural action of the group $\GL(m)$ on $\C_k$
defines canonically the action $\GL(m)_{k+2} : \C_k$
such that the subgroup $\N(m)_{k+2}$ acts trivially.
For $n = 0,1, \ldots , k$ let
\[
  \pr_n : \C_k \rightarrow C_n
\]
be canonical projections.
Then
\[
  \psi_n \overset{\text{def}}{ = }
  \pr_n \circ \widetilde{\psi}_k \ :
  \ \Con(m)_k \rightarrow C_n
\]
is $\GL(m)_{k+2}$-morphism.
From Theorem \ref{T:main2} it follows that $\psi_n$
defines canonically a natural differential operation
\[
  \Psi_n \ : \ \Con(M) \ \rightarrow
  \ \Gamma(\Ass_M(C_k)).
\]
From Lemma \ref{sesh2} we obtain the following statement.
\begin{theorem}
  Let $U$ be a $\GL(m)$-module and
  \[
    F \ : \ \Con(M) \ \rightarrow
    \ \Gamma(\Ass_M(U))
  \]
  be a natural differential operation of order $k$.
  Then
  \[
    F(\theta) = \widetilde{F}
    (\Psi_0(\theta), \ldots , \Psi_k(\theta)),
  \]
  where $\widetilde{F}$ corresponds to some polynomial
  $\GL(m)$-mapping
  \[
    \widetilde{f} \ : \ C_0 \times C_1 \times
    \ldots \times C_k \ \rightarrow \ U.
  \]
\end{theorem}
\begin{remark}
  From the construction of $\Psi_n$ we obtain the
  following geometrical description of the curvature
  $\Psi_n(\theta)|_p$ of the connection $\theta$ on $M$ at
  $p \in M$. Take local coordinates $(x^1, \ldots , x^m)$
  with center $p$ such that for the Taylor series
  \[
  \theta = \theta_0 + \theta_1 + \theta_2 + \ldots
  \]
  of the connection $\theta$ in the local coordinates
  $(x^1, \ldots , x^m)$ we have: $\theta_l \in C_l$
  for all $0 \le l \le n$. Then the local presentation
  at $p$ in the local coordinates $(x^1, \ldots , x^m)$ of
  the curvature $\Psi_n(\theta)$ is $\theta_n$.
\end{remark}
\begin{remark}
  It is easy to see that $\Psi_0(\theta)$ is the torsion of
  the connection $\theta$.
\end{remark}

\section{Natural differential operations on a manifold with
a connection.}\label{withconn}

\begin{theorem}\label{canonwithconn}
  Let $\rho : \GL(m) \rightarrow \GL(V)$ and
  $\sigma : \GL(m) \rightarrow \GL(U)$ be linear representations.
  Then there is the canonical 1-1 correspondence between
  $\GL(m)$-equivariant differential operations of bounded order
  in the form
  \begin{equation}\label{opwithconn}
    \xi : \ \C_k \times
    C^\infty(\R^m, V) \rightarrow C^\infty(\R^m, U),
  \end{equation}
  where $k \in \N$ and the set of natural differential operations
  of bounded order in the form
  \[
    \Xi : \Gamma(\Ass_M(V))
    \ \rightarrow \ \Gamma(\Ass_M(U))
  \]
  on a manifold with a connection.
\end{theorem}
\begin{proof}
  Let $z^1 , \ldots , z^m$ be the standard basis of the space
  $\R^{m\ast}$. We identify $\EE_k$ and the space of polynomials
  in the variables $z^1 , \ldots , z^m$ of degree $\le k$.
  
  Suppose we have $\GL(m)$-equivariant differential
  operations \eqref{opwithconn} of an order $\le k$.
  Consider the regular action
  \[
    \N(m)_{k+2} : \Con(m)_k \times (\EE_k \otimes V).
  \]
  Note that the subvariety $\C_k \times (\EE_k \otimes V)$
  is a nice Seshadri section for this action. Define
  $\GL(m)$-mapping
  \[
    \hat{\xi} \ : \ \C_k \times (\EE_k \otimes V)
    \rightarrow U, \qquad
    (\theta, v(z)) \mapsto \xi(\theta, v(z))(0).
  \]
  Extend the mapping $\hat{\xi}$ to the $\GL(m)_{k+2}$-mapping
  \[
    \xi \ : \ \Con(m)_k \times (\EE_k \otimes V)
    \ \rightarrow \ U,
  \]
  From Theorem \ref{T:main2} it follows that $\xi$ defines
  a natural differential operation
  \[
    \tilde{\Xi} \ : \ \Con(M) \times \Gamma(\Ass_M(V))
    \ \rightarrow \ \Gamma(\Ass_M(U)).
  \]
  Now suppose that $(M, \theta)$ is a manifold with
  a connection. Define the operation
  \[
    \Xi \ : \Gamma(\Ass_M(V))
    \ \rightarrow \ \Gamma(\Ass_M(U)), \qquad
    s \mapsto \tilde{\Xi}(\theta, s).
  \]

  Conversely, by the construction above any natural differential
  operation of bounded order on a manifold with a connection
  defines canonically $\GL(m)$-equivariant differential operation
  \eqref{opwithconn} of bounded order.
\end{proof}
\begin{remark}\label{toq}
  From the proof of Theorem \ref{mainonriem} we obtain
  the following description of the corresponding to $\xi$
  differential operation $\Xi$. Let $M$ be a manifold with
  a connection $\theta$, $s$ be a section of
  $\Ass_M(V)$, and $p \in M$. Take local coordinates
  $(x^1, \ldots , x^m)$ with center $p$ such that for the Taylor
  series $\theta = \theta_0 + \theta_1 + \ldots$
  of the connection $\theta$ in the local coordinates
  $(x^1, \ldots , x^m)$ we have: $\theta_l \in C_l$
  for all $0 \le l \le n$. Now calculate $\Xi(s)(p)$ in that
  local coordinates by the formula \eqref{opwithconn}.
\end{remark}

\section{Canonical quantization of the Poisson structures
on a manifold with a connection.}\label{q}

We need $\star(\omega)$-product in $\R^m$.

\begin{definition}
Let $x = (x^1, \ldots , x^m)$ be coordinate functions
in $\R^m$ and
\[
  \omega = \omega(x) = \omega^{ij}(x) \frac{\partial}{\partial x^i}
  \wedge \frac{\partial}{\partial x^j}
\]
be a Poisson structure. $\star(\omega)$-product in $\R^m$
is a mapping
\begin{gather*}
  \star(\omega) \ :
  \ C^\infty(\R^m) \times C^\infty(\R^m)
  \rightarrow C^\infty(\R^m)[[\hbar]], \\
  \star(\omega)(f,g) \ = \ f g + \beta_1(\omega,f,g) \hbar +
  \beta_2(\omega,f,g) \hbar^2 + \ldots
\end{gather*}
such that
\begin{itemize}
  \item $\star(\omega)$ is a quantization of the Poisson
    structure $\omega$,
  \item 
    \[
      \beta_k : C^\infty(\R^m, \wedge^2(\R^m)) \times
      C^\infty(\R^m) \times C^\infty(\R^m)
      \ \rightarrow \ C^\infty(\R^m)
    \]
is $\GL(m)$-equivariant differential operation of bounded order,
where $k = 1, 2, \ldots $.
\end{itemize}
\end{definition}

$\star(\omega)$-product in $\R^m$ is a natural generalization
of the Moyal $\star$-product. The first $\star(\omega)$-product
in $\R^m$ was constructed by Kontsevich in \cite{Ko}.
\begin{remark}
  By using the IT-reduction it is not difficult to prove that
  there exists $\star(\omega)$-product in $\R^m$.
  This proof is constructive: it gives an algorithm for
  calculation the operations $\beta_k$.
\end{remark}

\begin{theorem}
  For a manifold with a connection and a Poisson structure there
  exists the canonical quantization.
\end{theorem}
\begin{proof}
  Let $M$ be a manifold with a connection $\theta$ and
  a Poisson structure $\omega$. Let the natural differential
  operation
  \[
  B_k : \Gamma(\Ass_M(\wedge^2 \R^m)) \times
  C^\infty(M) \times C^\infty(M) \rightarrow C^\infty(M)
  \]
  corresponds to $\beta_k$, where $k = 1, 2, \ldots$
  (see Theorem \ref{canonwithconn}).
  Set
  \begin{equation}\label{form}\begin{split}
    \star \ : \ C^\infty(M) \times C^\infty(M)
    \rightarrow C^\infty(M)[[\hbar]], \\
    (f,g) \mapsto f g + B_1(\omega, f, g) \hbar +
    B_2(\omega, f, g) \hbar^2 + \ldots .
  \end{split}\end{equation}
  We claim that the operation \eqref{form} is a quantization
  of the Poisson structure $\omega$. To prove it we take
  a point $p \in M$ and $n \in \N$ and prove that the defined
  above operation $\star$ gives a quantization of the Poisson
  structure $\omega$ modulo $O(\hbar^{n+1})$. Suppose that
  $(x^1, \ldots , x^m)$ are local coordinates with center $p$
  such that for the Taylor series
  $\theta = \theta_0 + \theta_1 + \ldots$
  of the connection $\theta$ in the local coordinates
  $(x^1, \ldots , x^m)$ we have: $\theta_l \in C_l$
  for all $0 \le l \le n$. By Remark \ref{toq} the operation
  $\star$ at the point $p$ in the local coordinates
  $(x^1, \ldots , x^m)$ coincides with the Moyal-Kontsevich
  $\star(\omega)$-product modulo $O(\hbar^{n+1})$. This
  concludes the proof.
\end{proof}
\begin{remark}
  From section \ref{onriem} one can easily obtain the following
  quantization rule for a Riemannian manifold.
  
  \it Let $M$ be a Riemannian manifold with a Poisson
  structure $\omega$, $p \in M$. Take
  normal local coordinates $(x^1, \ldots , x^m)$ with
  center $p$. Then define $\star$-product at $p$ in coordinates
  $(x^1, \ldots , x^m)$ by the Moyal-Kontsevich formula for
  $\star(\omega)$-product. \rm
\end{remark}

\end{document}